\documentclass[12pt]{amsart}
\usepackage{graphicx}
\usepackage{amsmath}
\usepackage{amsthm,enumerate}
\newtheorem{thm}{Theorem}

\newtheorem{lem}[thm]{Lemma}
\newtheorem{prop}[thm]{Proposition}
\theoremstyle{definition}
\newtheorem{defn}[thm]{Definition}
\theoremstyle{remark}

\numberwithin{equation}{section}

\newcommand{\Rr}{\mathbb R}

\newcommand{\D}{\mathbb{D}}
\newcommand{\Hh}{\mathbb{H}}
\newcommand{\argg}{\operatorname{Arg\,}}
\begin{document}

\title[]{Interpolation by positive harmonic functions}

\author{D.\ Blasi and A.\ Nicolau}
\address{Daniel Blasi, Departament de Matem\`{a}tiques, Universitat Aut\`{o}noma de Barcelona, 08193 Bellaterra, Barcelona, Spain}
\urladdr{http://www.mat.uab.es/\textasciitilde dblasi/}
\email{dblasi@mat.uab.es}
\address{Artur Nicolau, Departament de Matem\`{a}tiques, Universitat Aut\`{o}noma de Barcelona, 08193 Bellaterra, Barcelona, Spain}
\urladdr{http://www.mat.uab.es/\textasciitilde artur/}
\email{artur@mat.uab.es}
\thanks{Both authors are supported in part by grants
MTM2005-00544 and 2001SGR00431\\
The first author is also supported by DURSI and Fons Social
Europeu under the grant 2003FI00116 }


\begin{abstract}
A natural interpolation problem in the cone of positive harmonic
functions is considered and the corresponding interpolating
sequences are geometrically described.
\end{abstract}

\maketitle

\section{Introduction}

Let $h^+=h^+ (\D)$ be the cone of positive harmonic functions in
the unit disc $\D$ of the complex plane. If $u\in h^+$ the
classical Harnack's inequality tells that
\begin{equation*}
\frac{1-|z|}{1+|z|} \leq \frac{u(z)}{u(0)} \leq
\frac{1+|z|}{1-|z|}
\end{equation*}
for any $z\in \D$. Recall that the hyperbolic distance
$\beta(z,w)$ between two points $z,w\in \D$ is
\begin{equation*}
\beta(z,w)=\log_2\,\frac{1+\left|\frac{z-w}{1-\bar w
z}\right|}{1-\left|\frac{z-w}{1- \bar w z}\right|}
\end{equation*}
Hence, estimates above can be read as $| \log_2 u(z)-\log_2 u(0)|
\!\leq\! \beta (z,0)$. Since these notions are preserved by
automorphisms of the disc, we deduce
\begin{equation}\label{1.1}
| \log_2 u(z)-\log_2 u(w)| \leq \beta (z,w)
\end{equation}
for any $z,w\in \D$. So for any function $u\in h^+$, a sequence of
points $\{z_n\}\subset \D$ and the  corresponding sequence  of
values $w_n=u(z_n)$, $n=1,2,\dotsb$ are linked by $|\log_2
w_n-\log_2 w_m|\leq \beta (z_n,z_m)$, $n,m=1,2,\dotsb$ However,
given a sequence of points $\{z_n\}\subset \D$, one can not expect
to interpolate by a function in $h^+$ any sequence of positive
values $\{w_n\}$ satisfying the above compatibility condition
unless the sequence $\{z_n\}$ reduces to two points. Actually it
is well known that having equality in (\ref{1.1}) for two distinct
points $z,w\in \D$ forces the function $u$ to be a Poisson kernel
and hence one can not expect to interpolate further values. In
other words, the natural trace space given by Harnack's
Lemma~(\ref{1.1}) is too large and one may consider the following
notion.

A sequence of points $\{z_n\}$ in the unit disc will be called an
interpolating sequence for $h^+$ if there exists a constant
$\varepsilon=\varepsilon (\{z_n\})>0$, 
such
that for any sequence of positive values $\{w_n\}$ satisfying
\begin{equation}\label{1.2}
 |\log_2 w_n-\log_2 w_m|\leq \varepsilon \beta (z_n,z_m), \quad n,m=1,2,\dotsb
\end{equation}
there exists a function $u\in h^+$ with $u(z_n)=w_n$, $n=1,2,\dotsb$

Observe that this is a conformally invariant notion, that is, if
$\{z_n\}$ is an interpolating sequence for $h^+$, so is $\{\tau
(z_n)\}$, for any automorphism $\tau$ of the unit disc. Moreover
the corresponding constants satisfy $\varepsilon (\{\tau
(z_n)\})=\varepsilon (\{z_n\})$. Recall that a sequence of points
$\{z_n\}$ in the unit disc is called separated if
$\inf\limits_{n\neq m} \beta (z_n,z_m)>0$. The main result of this
paper is the following.

\begin{thm}\label{thm1}
A separated sequence $\{z_n\}$ of points in the unit disc is
interpolating for $h^+$ if and only if there exist constants $M>0$
and $0<\alpha <1$ such that
\begin{equation}\label{1.3}
\# \{z_j \colon \beta (z_j, z_n)\leq l\}\leq M2^{\alpha l}
\end{equation}
for any $n,l=1,2,\dotsb$
\end{thm}

We have restricted attention to separated sequences because we
want to consider an interpolation problem by positive harmonic
functions and not by their derivatives. However it is worth
mentioning that any interpolating sequence for $h^+$  is the union
of at most three separated sequences. Let us now discuss
condition~(\ref{1.3}). As it is usual in this kind of problems,
the geometrical description of interpolating sequences is given in
terms of a density condition which tells, in the appropriate
sense, that interpolating sequences are not too dense. The number
~2 shows up in~ (\ref{1.3}) because of the normalization of the
hyperbolic distance. We have chosen this normalization because it
fits perfectly well with dyadic decompositions. As we will show in
Section~\ref{sectioncuatro}, there are a number of conditions
which are equivalent to ~(\ref{1.3}). For instance, a sequence
$\{z_n\}$ satisfies~(\ref{1.3}) if and only if there exist
constants $M_1>0$ and $0<\alpha <1$ such that
\begin{equation*}
\# \left \{ z_j \colon\left | \frac{z_j-z_n}{1-\bar z_nz_j}\right
|\leq r\right \} \leq M_1(1-r)^{-\alpha}
\end{equation*}
for any $n=1,2,\dotsb$ and $0<r<1$. One can also write an
equivalent condition in terms of Carleson measures. It will be
shown in Section~\ref{sectioncuatro} that a  sequence
$\{z_n\}\subset \D$ satisfies (\ref{1.3}) if and only if there
exist constants $M_2>0$ and $0<\alpha <1$ such that
\begin{equation*}
\sum_j (1-|z_j|)^\alpha \leq M_2(1-|z_n|)^\alpha , \quad n=1,2,\dotsb,
\end{equation*}
where the sum is taken over all points $z_j\in \{z_k\}$  such that
$|z_j-z_n|\leq 2(1-|z_n|)$. This resembles the usual Carleson
condition with an exponent $\alpha <1$ for the Carleson squares
which contain a point of the sequence in its top part. Let us now
discuss the geometrical meaning of condition~(\ref{1.3}). It tells
that, when viewed from a point of the sequence, sequences
satisfying (\ref{1.3}) are ---at large scales--- exponentially
more sparse than merely separated sequences. Actually, a sequence
of points $\{z_n\}\subset  \D$ is a finite union of separated
sequences if and only if (\ref{1.3}) holds with $\alpha =1$. It
should also be mentioned that in condition (\ref{1.3}) one counts
points in  the sequence which are at hyperbolic distance less than
$l$ from a given point $z_n$ in the sequence, instead of taking as
a base point any $z\in  \D$ as in \cite{BN}. See also [S,\ p.\
63--77]. This last condition is stronger. Actually it will be
shown in Section~\ref{sectioncuatro} that there exist two
separated interpolating sequences $Z_1,Z_2$ for $h^+$ with $\inf\{
\beta (z, \xi)\colon z\in Z_1, \xi\in Z_2\}>0$ such that $Z_1 \cup
Z_2$ is not an interpolating sequence for $h^+$.

Let us now explain the main ideas of the proof. Let $E^*$ denote the radial projection of a set $E\subset \D$, that is, $E^*=\{\xi \in \partial \D\colon r\xi \in E$ 
for some  $0\leq r<1\}$. An application of Hall's Lemma yields
that there exists a universal constant $C>0$ such that for any
$u\in h^+$ one has
\begin{equation*}
\left | \left \{ z\in \D\, \colon \frac{u(z)}{u(0)} >
\lambda\right \}^*\right | \leq \frac{C}{\lambda}, \quad \lambda
>0.
\end{equation*}
The necessity of condition (\ref{1.3}) follows easily from this
estimate. The proof of the sufficiency is harder. Given a sequence
of points $\{z_n\}\subset \D$ satisfying (\ref{1.3}) and a
sequence of positive values $\{w_n\}$ satisfying the compatibility
condition (\ref{1.2}), one has to find a function $u\in h^+$ such
that $u(z_n)=w_n$, $n=1,2,\dotsb$ The construction of the function
$u\in h^+$ may be splitted into three steps.

\medskip

1. We will apply a classical result in Convex Analysis called
Farkas Lemma which  may be understood as an analogue for Cones of
the Hahn-Banach Theorem. Instead of constructing directly the
function $u\in h^+$ which performs the interpolation, Farkas Lemma
will tell that it suffices to prove the following statement. Given
any partition of the sequence $\{z_n\}$ into two disjoint
subsequences, $\{z_n\}=T\cup S$, there exists a function
$u=u(T,S)\in h^+$ such that
\begin{equation*}
\begin{split}
u(z_n)\geq w_n, \quad \textrm{ if } z_n\in T,\\
u(z_n)\leq w_n, \quad \textrm{ if } z_n\in S.
\end{split}
\end{equation*}

\medskip

2. Let $\omega(z,G)$ denote the harmonic measure in $\D$ of the
set $G\subset \partial \D$ from the point $z\in \D$, that is,
\begin{equation*}
\omega(z,G)=\frac{1}{2\pi} \int_G \frac{1-|z|^2}{|\xi-z|^2} \,
|d\xi|.
\end{equation*}
For each point $z_n$ of the sequence $\{ z_n\}$ we will construct
a set $G_n$ $\subset
\partial \D$ and  we will show that condition (\ref{1.3}) provides
some sort of independence of harmonic measures $\{\omega(z_n,
\cdot)\colon n=1,2,\dotsb\}$. Actually, given $0<\delta <1$, there
exists $N>0$ and a collection of pairwise disjoint subsets
$\{G_n\}$ of $\partial \D$ such that
\begin{equation*}
\begin{split}
&\omega(z_n, \cup_{k\in A(n)} G_k) \geq 1-\delta, \\
&\sum_{k\notin A(n)} 2^{\eta \beta (z_k, z_n)} \omega(z_n,
G_k)\leq \delta.
\end{split}
\end{equation*}
Here $A(n)$ denotes the set of indexes $k$ so that $\beta (z_k,
z_n)\leq N$. The number $\eta=\eta(\delta, M, \alpha)>0$ is a
constant depending on $\delta >0$ and on the constants $M>0$ and
$\alpha <1$ of (\ref{1.3}). The construction of the sets $\{G_n\}$
uses a certain stopping time argument and constitutes the most
technical part of the proof.

\medskip

3. L.\ Carleson ad J.\ Garnett found a description of the interpolating
sequences for the space $h^\infty$ of bounded harmonic functions in the unit
  disc (see \cite{CG}, \cite{G1} or \cite[p.\ 313]{G2}). Using their result it is easy to show that a separated sequence
  verifying (\ref{1.3}) is interpolating for $h^\infty$. Hence there exists
   $\gamma >0$ and a harmonic function $h$, with $\sup\{|h(z)|
   \colon z\in \D\} <1$ such that $h(z_n)=\gamma$ 
if $z_n\in T$, while $h(z_n)=-\gamma$ 
if $z_n\in S$. Then, fixed $\varepsilon >0$ and $\delta >0$
sufficiently small, using the compatibility condition (\ref{1.2})
and the estimates in step 2, one can show that the function
\begin{equation*}
u(z)=\sum_{z_n\in T} w_n \int_{G_n} \frac{1-|z|^2}{|\xi-z|^2}\, (1+h(\xi))
\frac{|d\xi|}{2\pi}, \quad z\in \D,
\end{equation*}
verifies $u(z_n)\geq w_n$ if $z_n\in T$ and $u(z_n)\leq w_n$ if $z_n\in S$.

\medskip

One may consider a similar problem in higher dimensions. Let
$h^+(\Rr_+^{d+1})$ denote the cone of positive harmonic functions
in the upper half space $\Rr_+^{d+1}=\{(x,y)\colon x\in \Rr^d,
y>0\}$. A sequence of points $\{z_n\}\subset \Rr_+^{d+1}$ will be
called an interpolating sequence for $h^+(\Rr_+^{d+1})$ if there
exists a constant $\varepsilon=\varepsilon(\{z_n\}) >0$ such that
for any sequence of positive values $\{w_n\}$ verifying
\begin{equation*}
 |\log_2 w_n-\log_2 w_m|\leq \varepsilon \beta (z_n,z_m), \quad n,m=1,2,\dotsb
\end{equation*}
there exists $u\in h^+(\Rr_+^{d+1})$ with $u(z_n)=w_n$,
$n=1,2,\dotsb$. When $d>~1$ we do not have a complete geometric
description of interpolating sequences. In this direction the
situation is analogue to the work of L.~Carleson and J.\ Garnett
\cite{CG} on interpolating sequences for the space $h^\infty
(\Rr_+^{d+1})$ of bounded harmonic functions in $\Rr_+^{d+1}$. See
section~\ref{sectionseis} for details.

The paper is organized as follows: Section~\ref{sectiondos} is
devoted to the proof of the necessity of condition (\ref{1.3}).
Section~\ref{sectiontres} contains the proof of the sufficiency.
Section \ref{sectioncuatro} is devoted to the analysis of
condition (\ref{1.3}). In Section ~\ref{sectioncinco} a related
interpolation problem for bounded analytic functions in the unit
disc without zeros is considered. This may be compared to
\cite{DN}. In the last section the interpolation problem for
positive harmonic functions in higher dimensions is discussed. The
letter $C$ will denote an absolute constant whose value may change
from line to line. Also $C(M)$ will denote a constant which
depends on $M$.

\section{Necessity}\label{sectiondos}

Given a set $E\subset \D$, let $\omega (z,E,\D\setminus E)$ denote
the harmonic measure from the point $z\in \D\setminus E$ of the
set $E$ in the domain $\D\setminus E$. The classical Hall's Lemma
tells that there exists a universal constant $C>0$ such that
$\omega (0,E,\D\setminus E)\geq C|E^*|$ for any set $E\subset \D$.
See \cite{H} or \cite{MS}. Recall that $E^*$ denotes the radial
projection of $E$. The main auxiliary result is the following.

\begin{lem}\label{lemma21}
There exists a constant $C>0$ such that for any $u\in h^+$ and $\lambda >0$ one has
\begin{equation*}
\left | \left \{ z\in \D\colon \frac{u(z)}{u(0)} > \lambda\right
\}^*\right | \leq \frac{C}{\lambda}
\end{equation*}
\end{lem}

\begin{proof}
One may assume that $\lambda >1$. Fix $u\in h^+$, let $E=\{z\in
\D\colon$ \linebreak $u(z)>\lambda u(0)\}$. The maximum principle
shows that
\begin{equation*}
u(z)\geq \lambda u(0)\omega (z,E,\D\setminus E), \quad z\in
\D\setminus E.
\end{equation*}
Taking $z=0$, one gets $\omega(0,E,\D\setminus E)\leq
\lambda^{-1}$ and applying Hall's Lemma one finishes the proof.
\end{proof}

\begin{proof}[Proof of the necessity of condition (\ref{1.3})]
Assume that $\{z_k\}$ is an interpolating sequence for $h^+$. By
conformal invariance it is sufficient to prove  (\ref{1.3})   when
the base point $z_n$ is the origin. So assume $z_1=0$ and take
$w_k=2^{\varepsilon \beta (z_k,0)}$, $k=1,2,\dotsb$ It is clear
that the compatibility condition  (\ref{1.2}) holds. So, there
exists $u\in h^+$ with $u(z_k)=w_k$, $k=1,2,\dotsb$ Let $D_k$ be
the hyperbolic disc centered at $z_k$ of hyperbolic radius $1$. By
Harnack's Lemma
 \begin{equation*}
u(z)\geq \frac{w_k}{2}, \quad z\in D_k, \quad k=1,2,\dotsb
\end{equation*}
So, if $A(j)$ denotes the set of indexes $k$ corresponding to
points $z_k$ with $j-1\leq \beta (z_k,0)\leq j$, $j=1,2,\dotsb$,
one deduces
 \begin{equation*}
u(z)\geq 2^{\varepsilon (j-1)-1}, \quad z\in D_k, \quad k\in A(j).
\end{equation*}
Now since $u(0)=1$, Lemma \ref{lemma21} gives
\begin{equation*}
\left |\left( \cup_{k\in A(j)} D_k\right)^*\right | \leq C_1
2^{\varepsilon (1-j)}.
\end{equation*}
Since the sequence $\{z_k\}$ is separated, the discs $\{D_k\}$ are quasidisjoint and one deduces
\begin{equation*}
\sum_{k\in A(j)} 1-|z_k| \leq C_2 2^{\varepsilon (1-j)}.
\end{equation*}
Since $1-|z_k|$ is comparable to $2^{-j}$ for any $k\in A(j)$, one deduces
\begin{equation*}
\# A(j)\leq C_3 2^{(1-\varepsilon)j}.
\end{equation*}
Adding up for $j=1,\dotsb,l$, one gets
\begin{equation*}
\# \{z_k \colon \beta (z_k,0)\leq l \}\leq C_4
2^{(1-\varepsilon)l}. \;\qed
\end{equation*}  \def\qed{}
\end{proof}

\section{Sufficiency of Condition (\ref{1.3})}\label{sectiontres}

By a normal families argument, one may assume the sequence
$\{z_n\}$ consists of finitely many points. As explained in the
introduction the proof consists of three steps.

\subsection{First Step}\label{3.1} Let $e_1,\dotsb,e_m$ be a
collection of vectors of the
 euclidian space $\Rr^d$. Farkas Lemma asserts that a vector $e\in\Rr^d$
  is in the cone generated by $\{ e_i\colon i=1,\dotsb, m\}$, that is
  $e=\sum \lambda_i e_i$ for some $\lambda_i\geq 0$, $\; i=1,\dotsb,m$,
   if and only if $\langle x,e\rangle\leq 0$ for any vector $x\in \Rr^d$
    for which $\langle x,e_i\rangle\leq 0$, $i=1,\dotsb,m$. See \cite{HL}.
     This classical result will be used in the proof of the next auxiliary
      result

\begin{lem}\label{lemma31}
Let $\{z_n\}$  be a sequence of distinct points in the
unit disc and let $\{w_n\}$ be a sequence of positive values.
Assume that for every partition of the sequence $\{z_n\}=T \cup S$, into
two disjoint
 subsequences $T$ and $S$, there exists $u=u(T,S)\in h^+$ such that
 $u(z_n)\geq w_n$ if $z_n \in T$ and $u(z_n)\leq w_n$ if $z_n\in S$.
 Then, there exists $u\in h^+$ such that $u(z_n)=w_n$,  $n=1,2,\dotsb$
\end{lem}

\begin{proof}[Proof of Lemma \ref{lemma31}]
By a  normal families argument, one may assume that both the sequences
 of points $\{z_n\}$ and values $\{w_n\}$ consist of finitely many, say $d$,
  points.
Consider the set of  all  partitions $\{z_n\}=T_k \cup S_k$,
 $k=1,\dots,m$ of the sequence $\{z_n\}$. Let $u_1,\dots,u_m \in h^+$ be the
  corresponding functions such that
 $u_k(z_n)\geq w_n$ if $z_n\in T_k$ and $u_k(z_n)\leq w_n$ if $z_n\in S_k$,
  and consider the vector
\begin{equation*}
u_i \colon =(u_i(z_1),\dots, u_i(z_d)), \quad i=1,\dots,m.
\end{equation*}
If $x\in \Rr^d$ satisfies $\langle x, u_i \rangle\leq 0$,
$i=1,\dots,m$, that is $\sum_{n=1}^d u_i(z_n)x_n\leq 0$, let
$\mathcal{F}=\{z_n\colon x_n\geq 0\}$. Then $\mathcal{F}=T_k$ for
some $k$ and let $S_k=\{z_n\}\setminus \mathcal F$. Its
corresponding function $u_k$ satisfies $x_n w_n\leq x_n u_k(z_n)$
for all $n=1,\dotsb,d$. So,
\begin{equation*}
\langle x,w\rangle =\sum_{n=1}^d w_nx_n\leq \sum_{n=1}^d
u_k(z_n)w_n\leq 0.
\end{equation*}
Now, by Farkas's Lemma, $w=(w_1,\dotsb,w_d)$ is in the cone
generated by the vectors $\{u_i, i=1,\dots,m\}$. So there exist
constants $\lambda_i\geq 0$, $\; i=1,\dotsb, m$ such that
$u(z)=\sum_{i=1}^m\lambda_i u_i(z)\in h^+$ and $u(z_n)=w_n$,
$n=1,2,\dots,d. $
\end{proof}

\subsection{Second Step}\label{3.2} The second step in the proof
consists on using condition (\ref{1.3}) to construct a collection
of disjoint subsets $\{G_n\}$ of the unit circle which provide a
suitable kind of independence of the harmonic measures
$\{\omega(z_n, \cdot) \colon n=1,2,\dotsb\}$. The precise
statement is given in the following result which is the main
technical part of the proof. Recall that $\omega(z,G)$ denote the
harmonic measure in $\D$ of the set $G\subset \partial \D$ from
the point $z\in \D$, that is,
\begin{equation*}
\omega(z,G)=\frac{1}{2\pi} \int_G \frac{1-|z|^2}{|\xi-z|^2} \,
|d\xi|.
\end{equation*}

\begin{lem}\label{lemma33}
Let $\{z_n\}$ be a sequence of distinct points in the unit disc
which satisfies condition (\ref{1.3}). Then for any $\delta>0$,
there exist numbers $N=N(\delta)>0$, $\eta=\eta (\delta)>0$ and a
collection $\{G_n\}$ of pairwise disjoint subsets of the unit
circle such that
\begin{equation}\label{3.1}
\omega\left (z_n, \cup_{k\in A(n)}G_k\right )\geq 1-\delta, \quad
n=1,2,\dotsb,
\end{equation}
and
\begin{equation}\label{3.2}
\sum_{k\notin A(n)} 2^{\eta\beta(z_k,z_n)} \omega(z_n,G_k)<
\delta, \quad n=1,2,\dotsb.
\end{equation}
Here $A(n)=A(n,N)$ denotes the collection of indexes $k$ such that
$\beta (z_k,z_n)\leq N$.
\end{lem}

We first introduce some notation. Given  a point  $z\in \D$ and
$C>0$ we denote
\begin{equation*}
\begin{split}
I(z)  &=\{e^{i\theta}    \colon -\pi (1-|z|)<
\theta-\operatorname{Arg} z \leq \pi (1-|z|)\},\\*[5pt]
Q(z)&=\left \{ re^{i\theta}    \colon  0<1-r\leq 1-|z|,\;
e^{i\theta} \in I(z)\right\},\\
CI(z)&= \left\{ e^{i \theta} \colon -\pi \, C (1-|z|)<\theta -
\operatorname{Arg} z \leq \pi \, C(1-|z|)\right \}\\*[5pt]
CQ(z)&=\left \{ re^{i\theta}    \colon  0<1-r\leq C( 1-|z|),\;
e^{i\,\theta}\in C\; I(z)\right \}
\end{split}
\end{equation*}

Observe that if $C(1-|z|)\geq 1$, one has $CI(z)=\partial \, \D$
and $CQ(z)=\D$. When $z=z_k\in \{z_n\}$, we simply denote
$I_k=I(z_k)$. We will use the following two elementary auxiliary
results.

\begin{lem}\label{lemma34}
Fixed $\delta >0$, there exists $M_0=M_0(\delta)>0$  such that
\begin{equation*}
\omega(z_k, M_0I_k)\geq 1-\frac{\delta}{100}, \quad k=1,2,\dotsb
\end{equation*}
\end{lem}

\begin{proof}
If $z_k=0$ one may take $M_0=1$. If $z_k\neq 0$ observe that there
exists an absolute constant $C_0>0$ such that $|e^{it}-z_k|\geq
C_0|t-\operatorname{Arg} z_k|$. Since
\begin{equation*}
\omega(z_k,\partial \D \setminus
M_0I_k)=\frac{1-|z_k|^2}{2\pi}\int_{\partial \D\setminus M_0I_k}
\frac{|d\xi|}{|\xi-z_k|^2},
\end{equation*}
one gets
\begin{equation*}
\omega(z_k,\partial \D \setminus M_0I_k)\leq \frac{1-|z_k|^2}{2\pi
C_0^2}\int_{\pi M_0(1-|z_k|)}^\infty \frac{dx}{x^2}.
\end{equation*}
Hence
\begin{equation*}
\omega(z_k,\partial \D \setminus M_0I_k)\leq\frac{1}{\pi^2C_0^2
M_0}
\end{equation*}
and taking $M_0=100/\pi C_0^2 \delta$ the result follows.
\end{proof}

\begin{lem}\label{lemma35}
Fixed $M>0$, there exists a constant $C(M)>0$ such that for all
pair of points $z,w\in \D$ with  $w\in 20 MQ(z)$, one has
\begin{equation*}
\left |\beta (z,w)-\log_2 \left ( \frac{1-|z|}{1-|w|}\right ) \right | \leq C(M).
\end{equation*}
\end{lem}

\begin{proof}
One may assume that $z,w\in \D\setminus \{0\}$. Since
\begin{equation*}
|1-\bar w z| \geq \left (1-|z||w|\right ) \geq \left (1-|z|\right
)
\end{equation*}
and
\begin{equation*}
\begin{split}
|1-\bar w z|&\leq |w|\left | \frac{1}{\bar w}-z\right |\\
&\leq |w| \left | \frac{1}{\bar w}-e^{i\operatorname{Arg} w}\right
|+\left | e^{i\operatorname{Arg} w}- e^{i\operatorname{Arg}
z}\right | +
\left | e^{i\operatorname{Arg} z}- z\right | \\
&\leq (20M+20M\pi +1)(1-|z|),
 \end{split}
\end{equation*}
we deduce
\begin{equation*}
1-|z|\leq |1-\bar wz|\leq K(M)(1-|z|),
\end{equation*}
where $K(M)=20M+20M\pi+1$. So,
\begin{equation*}
\begin{split}
\beta (z,w) &=2\log_2\left (1+\left | \frac{z-w}{1-\bar wz}\right|\right )
-\log_2\left (1-\left | \frac{z-w}{1-\bar wz}\right|^2\right )\\
&= 2\log_2 \left ( 1+\left | \frac{z-w}{1-\bar wz}\right | \right ) -
\log_2  \frac{(1-|z|^2)(1-|w|^2)}{|1-\bar wz|^2} \\
&= C+\log_2
\left ( \frac{1-|z|}{1-|w|}\right)
 \end{split}
\end{equation*}
where $-2\leq C\leq 2+2\log_2 K(M)$.
\end{proof}



\begin{proof}[Proof of Lemma \ref{lemma33}]
The construction of the sets $\{G_n\}$ may be splitted into three
steps.
\begin{enumerate}[i)]
\item For each $z_k\in \{z_n\}$ and $\lambda >0$, we will
construct certain points $z_n^\gamma (k)\in \D$ with
$I(z_n)\subset I(z_n^{(\gamma)}(k))$ and
\begin{equation}\label{6.1}
\sum_{{\small \begin{array}{l} z_n\in 20M_0Q(z_k)\\
\beta(z_k,z_n)\geq N \end{array}}} 1-|z_n^\gamma (k)|\leq \lambda
(1-|z_k|) \textrm{ for all } z_k \in \{z_n\}.
\end{equation}
Here $N$ is a constant depending on $\lambda, M_0$ and on the
constants $M$ and $\alpha$ appearing in (\ref{1.3})

\item Next, we will construct certain
 sets $E_k\subset \partial \D$ with
$E_k\cap E_j=\emptyset$ if $\beta (z_k, z_j)\geq N$ such that
\begin{equation}\label{6.2}
\omega(z_k,E_k)\geq 1-\frac{\delta}{10}.
\end{equation}
In the construction of the sets $E_k$ we will use the points
$z_n^\gamma (k)$ of the first step which satisfy the estimate
(\ref{6.1}) above for a certain fixed $\lambda$ sufficiently
small.

\item Finally we will construct the pairwise disjoint sets $G_n$ satisfying
conditions (\ref{3.1}) and (\ref{3.2}).
\end{enumerate}

\bigskip

{\bf i) Construction of the points \boldmath{$z_n^\gamma (k)$}.} 
Fix $\delta >0$. Applying Lemma~\ref{lemma34}, there exists a
constant $M_0=M_0(\delta)>0$ such that
\begin{equation}\label{6.3}
\omega(z_k,M_0I_k)\geq 1-\frac{\delta}{100}, \quad k=1,2,\dotsb
\end{equation}
Fix $z_k\in \{z_n\}$. Let $\gamma=\gamma(\alpha)>0$ be a small
number to be fixed later. For any $z_n\in 20M_0Q(z_k)$ with $\beta
(z_k,z_n)\geq~N$ we define $z_n^\gamma (k)$ as the point in $\D$
satisfying the following three conditions
 \begin{equation}\label{6.4}
\begin{split}
& \operatorname{Arg} (z_n) =\operatorname{Arg} ( z_n^\gamma(k) ),\\
& \beta (z_n^\gamma (k),z_n)=\gamma\beta (z_k,z_n),\\
& \left | z_n^\gamma (k)\right | <|z_n|.
 \end{split}
\end{equation}
Here $N=N(\gamma, M_0,\lambda)$ is a large number to be fixed
later. In particular $N>0$ will be taken so large that
$z_n^{\gamma}(k)\in 20M_0Q(z_k)$ whenever $z_n\in 20M_0Q(z_k)$
satisfies $\beta (z_n,z_k)>N$. See Figure 1.

\begin{figure}[h]
\scalebox{0.75}{\includegraphics{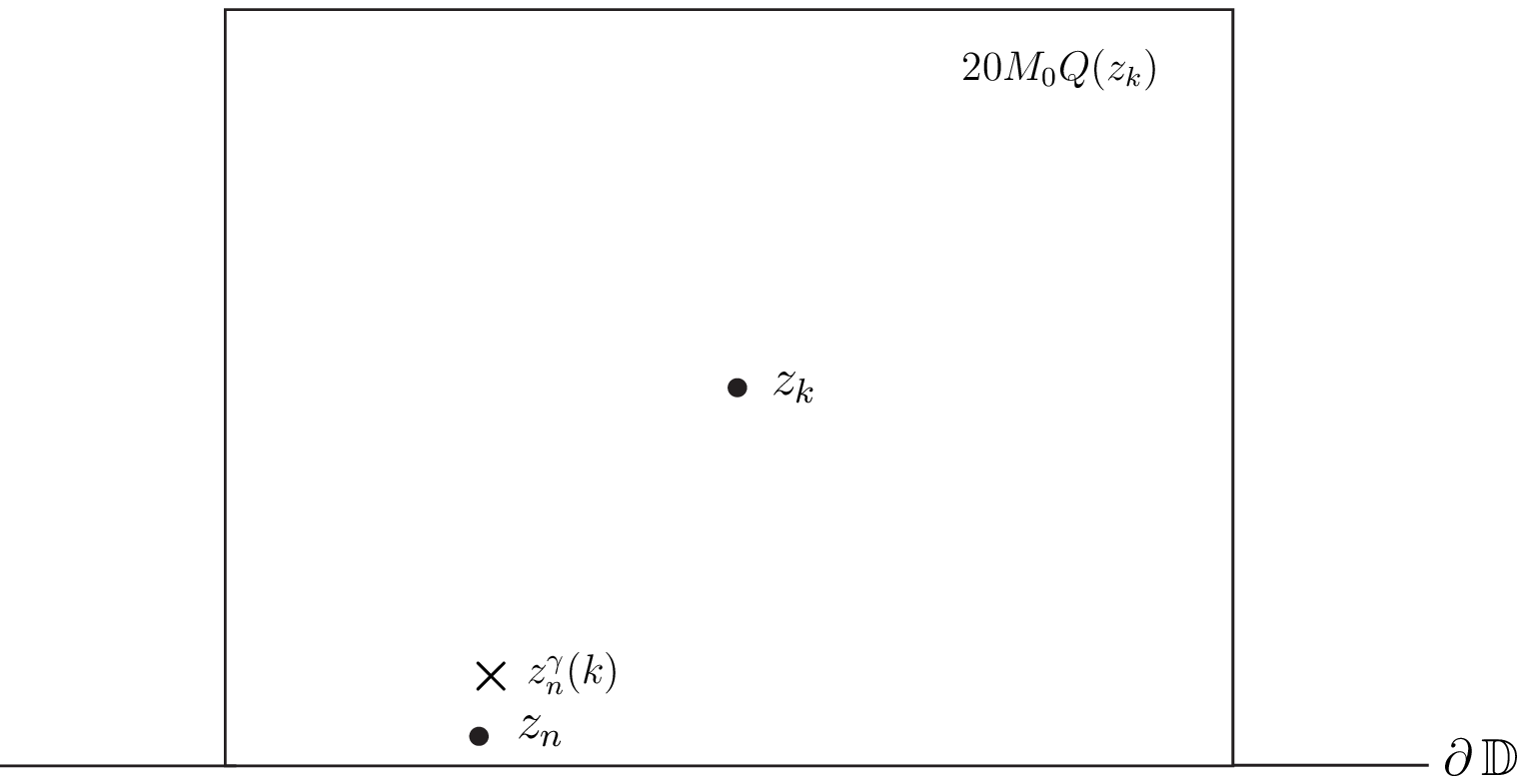}} \caption{}
\end{figure}

Using Lemma \ref{lemma35} and
$\beta(z_n^{\gamma}(k),z_n)=\gamma\beta (z_k,z_n)$
 we obtain the following inequalities:
 \begin{equation}\label{6.5}
\left (\frac{1-|z_k|}{1-|z_n|}\right )^{C^{-1}\gamma}\leq
\frac{1-|z_n^\gamma (k)|}{1-|z_n|}
 \leq \left ( \frac{1-|z_k|}{1-|z_n|}\right )^{C \gamma},
\end{equation}
where $C$ is a constant depending on $M_0$. So,
\begin{equation*}
\sum_{\begin{array}{c}
\small{
z_n\in 20M_0Q(z_k)}\\  \small{ \beta(z_k,z_n)\geq N} \end{array}} \hspace{-1cm}
 1-|z_n^\gamma (k)|\leq
\left (1-|z_k|\right )^{C\gamma} \sum_{j=N}^\infty \hspace{-.5cm}
\sum_{{\small \begin{array}{c}
z_n\in 20M_0Q(z_k)\\ j\leq \beta(z_n,z_k)<j+1 \end{array}}} \hspace{-1cm}
\left (1-|z_n|\right)^{1-C\gamma}.
\end{equation*}
Now, if $z_n\in 20 M_0Q(z_k)$ and $j\leq \beta (z_n,z_k)< j+1$,
Lemma~\ref{lemma35} tells that
 $1-|z_n|\leq K(M_0)2^{-j}(1-|z_k|)$. So, using 
(\ref{1.3}), the right hand side term is bounded by
\begin{equation*}
K(M_0)^{1-C\gamma}(1-|z_k|)\sum_{j=N}^\infty M2^{\alpha j}
2^{-j(1-C\gamma)}.
\end{equation*}
Since $\alpha <1$, taking $\gamma >0$ so small that
$\alpha+C\gamma<1$, the expression above may be bounded by
\begin{equation*}
M\, K(M_0)^{1-C\gamma}\frac{2^{N(\alpha
+C\gamma-1)}}{1-2^{\alpha+C\gamma-1}}\, (1-|z_k|).
\end{equation*}
Finally, given $\lambda>0$ taking $N$ sufficiently large, we obtain
\begin{equation*}
\sum_{\begin{array}{c}
\small{
z_n\in 20M_0Q(z_k)}\\  \small{ \beta(z_n,z_k)\geq N} \end{array}}  1-|z_n^\gamma (k)|\leq \lambda \left (1-|z_k|\right ) \textrm{ for all } z_k \in \{z_n\}.
\end{equation*}

\bigskip

{\bf ii) Construction of the sets \boldmath{$\{E_k\}$}.} For each
$z_n^\gamma(k)$, we define $I_n^\gamma (k)=I(z_n^{(\gamma)}(k))$.
Fixed $M_0>0$ and $N>0$, we introduce the notation:
\begin{equation*}
B(k)= \{ z_n \;\colon \; |z_n|\geq |z_k|, \; \beta (z_k,z_n)\geq N, \; z_n\in 20M_0Q(z_k)  \}.
\end{equation*}
Now we will proof that the sets $E_k=M_0I_k\setminus
\bigcup_{z_n\in B(k)} I_n^\gamma (k)$ satisfy
\begin{equation}\label{eq6x}
\omega(z_k,E_k)\geq
1-\frac{\delta}{10}.
\end{equation}
Using the elementary estimate of the Poisson Kernel
\begin{equation*}
\frac{1-|z_k|^2}{|e^{it}-z_k|^2} \leq \frac{1+|z_k|}{1-|z_k|},
\end{equation*}
one obtains
\begin{equation*}
\omega\bigl(z_k, \hspace*{-.2cm}\bigcup_{z_n\in B(k)}\! I_n^\gamma (k)\bigr ) 
\!\leq \hspace*{-.2cm}
\sum_{z_n\in B(k)} \frac{1+|z_k|}{1-|z_k|}
\int_{I_n^\gamma (k)} \frac{dt}{2\pi}\\
 \leq \frac{2}{1-|z_k|}
\sum_{z_n \in B(k)} \!\!1-|z_n^\gamma(k)|.
\end{equation*}
which by (\ref{6.1}) is smaller than $2\lambda$. Since
\begin{equation*}
\omega(z_k,E_k)=\omega(z_k,M_0I_k)-\omega\left (z_k,
\bigcup_{z_n\in B(k)} I_n^\gamma(k)\right ),
\end{equation*}
the estimate (\ref{6.3}) tells
\begin{equation*}
\omega(z_k,E_k)\geq 1-\frac{\delta}{100}-\lambda.
\end{equation*}
If we take $\lambda>0$ sufficiently small, we deduce (\ref{eq6x}).
Since $M_0 I_n\subset I_n^\gamma(k)$, it is clear from the
definition that $E_k\cap E_j=\emptyset$ if $\beta (z_k,z_j)>N$.

\medskip

{\bf iii) Construction of the pairwise disjoint sets
\boldmath{$G_n$}.} We rearrange the sequence $\{z_n\}$ so that
$\{1-|z_n|\}$ decreases. For each point $z_n$ we will construct a
set $G_n\subset E_n$ so that the corresponding family $\{G_n\}$
will satisfy (\ref{3.1}), (\ref{3.2})  and $G_n \cap
G_m=\emptyset$ if $n\neq m$. The construction will proceed by
induction and will ensure that the sets $G_n$ are pairwise
disjoint and verify (\ref{3.1}).

Take $G_1=E_1$. By (\ref{eq6x}), the estimate (\ref{3.1})
is satisfied when $n=1$. Assume that pairwise disjoint subsets
$G_1,\dotsb, G_{j-1}$ of the unit circle have been defined so that
\begin{equation*}
\omega (z_n,\hspace*{-.2cm} \bigcup_{k\leq n,k\in A(n)} \!G_k)\geq 1-\delta, \textrm{ for } n=1,2,\dotsb, j-1.
\end{equation*}
The set $G_j$ will be constructed according to the following two different situations:

\begin{enumerate}[(1)]
\item If $\beta (z_j, \{z_1, \dotsb, z_{j-1}\}) \geq N$ we define
$G_j=E_j$. By (\ref{6.2}) we have
\begin{equation*}
\omega (z_j, \bigcup_{k\leq j,k\in A(j)} G_k)\geq
\omega(z_j,G_j)\geq 1-\delta.
\end{equation*}
Now let us show that $G_k\cap G_j=\emptyset$ for any $k=1,\dotsb,j-1$. Since $G_k\subset E_k$ and $G_j \subset M_0 I_j$, it is sufficient to show that
$M_0 I_j\cap E_k=\emptyset$ for $k=1,\dotsb,j-1$. Fix $k=1,\dotsb, j-1$ and consider two cases
\begin{enumerate}[(a)]
\item If $z_j\in 20 M_0Q(z_k)$. Since $M_0 I_j\subset I_j^\gamma (k)$ and $E_k=M_0 I_k\setminus \bigcup I_j^\gamma(k)$, we have $E_k\cap M_0 I_j=\emptyset$
\item If $z_j\notin 20 M_0 Q(z_k)$. Since $|z_j|>|z_k|$ we have $M_0 I_j \cap M_0 I_k =\emptyset$. Hence $E_k\cap M_0 I_j=\emptyset$.
\end{enumerate}

\item If $\beta (z_j, \{z_1,\dotsb, z_{j-1}\})\leq N$,
consider the set  of indexes $\mathcal{F}=\mathcal F(j)=\{ k\in
[1,\dotsb, j-1] \colon \beta (z_k,z_j)\leq N\}$. Let us
distinguish the following two cases:
\begin{enumerate}[(a)]
\item If $\omega(z_j, \bigcup_{k\in \mathcal{F}} G_k)\geq
1-\delta$, define $G_j=\emptyset$. It is obvious that
\begin{equation*}
\omega (z_j, \bigcup_{k\leq j, k\in A(j)} G_k)\geq 1-\delta.
\end{equation*}
\item If $\omega(z_j, \bigcup_{k\in \mathcal{F}} G_k)< 1-\delta$,
define $G_j=E_j\setminus \, \bigcup_{k\in \mathcal{F}} G_k$.
Arguing as in case 1 one can show that $G_k\cap G_j=\emptyset$ for any $k=1,\dotsb,j-1$. Also, applying (\ref{eq6x}),  one gets
\begin{equation*}
\omega (z_j, \bigcup_{k\leq j, k\in A(j)} G_k)\geq  \omega (z_j,E_j) \geq1-\delta.
\end{equation*}
\end{enumerate}
\end{enumerate}
So, by induction, a family $\{G_n\}$ of pairwise disjoint subsets
of the unit circle is constructed so that condition (\ref{3.1}) is
satisfied. It just remains to show that the family $\{G_n\}$
verifies (\ref{3.2}), that is, there exists $\eta=\eta (\delta)>0$
such that
\begin{equation*}
\sum_{k\colon \beta (z_k,z_n)\geq N} 2^{\eta \beta (z_k,z_n)}
\omega(z_n,G_k)\leq \delta, \quad n=1,2,\dotsb
\end{equation*}
Fixed $ n=1,2,\dotsb$, split this sum into three parts (A), (B)
and (C), corresponding to the points $z_k$ with $\beta
(z_k,z_n)\geq N$ such that:
\begin{enumerate}[(a)]
   \item $z_k\in 20 M_0 Q(z_n)$ in part (A),
   \item $z_k$ so that $2M_0 I_k\cap M_0 I_n=\emptyset$ in part   (B)\qquad (See Figure 2)
   \item points $z_k$ which are not in (a) or (b)
\end{enumerate}

\begin{figure}[h]
\scalebox{0.75}{\includegraphics{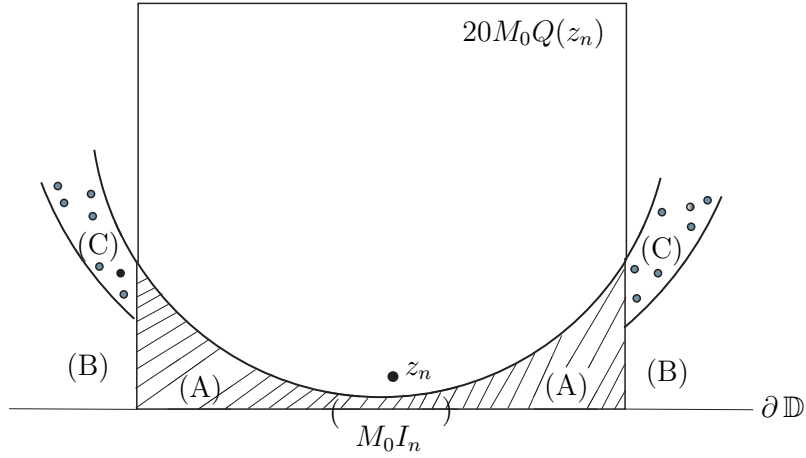}} \caption{The sum is
splitted into three parts corresponding to the location of the
points $z_k$ in the regions $\operatorname{(A)}$,
$\operatorname{(B)}$ or  $\operatorname{(C)}$}
\end{figure}


In (A) and (B) we will use the estimate $\omega(z_n,G_k)\leq
C(M_0)2^{-\beta (z_n,z_k)}$ and for (C) we will use the constant
$\gamma
>0$ appearing in the construction of the sets $E_k$.

We first claim that there exists a constant $C= C(M_0)>0$ such
that for points $z_k$ in part (A) or (B), that is those verifying
either  $z_k\in 20 M_0 Q(z_n)$ or $2M_0 I_k\cap M_0
I_n=\emptyset$, one has
\begin{equation}\label{eq3.10}
\omega (z_n,G_k)\leq C 2^{-\beta(z_k,z_n)}.
\end{equation}
For the points $z_k$ in part (A) we have $z_k\in 20 M_0 Q(z_n)$. Since $G_k\subseteq M_0 I_k$, a trivial estimate of the Poisson kernel gives
\begin{equation*}
\omega (z_n,G_k)\leq \int_{M_0 I_k} \frac{1-|z_n|^2}{|e^{it}-z_n|^2} \; \frac{dt}{2\pi} \leq 2M_0 \frac{1-|z_k|}{1-|z_n|}.
\end{equation*}
Applying Lemma \ref{lemma35}, since $z_k\in 20 M_0Q(z_n)$, one has
\begin{equation*}
\log_2 \frac{1-|z_k|}{1-|z_n|}\leq C (M_0)-\beta (z_k,z_n).
\end{equation*}
Hence, if $z_k\in 20 M_0 Q(z_n)$ we deduce
\begin{equation*}
\omega (z_n,G_k)\leq C2^{-\beta (z_k,z_n)}
\end{equation*}
with $C=2M_0 2^{C(M_0)}$. For the points $z_k$ in part (B) we have $2M_0 I_k\cap M_0 I_n=\emptyset$. An easy calculation shows that there exists a constant $C_1=C_1(M_0)$ such that for any $e^{it} \in I_k$ one has
\begin{equation*}
|e^{it} -z_n|\geq C_1 |1-z_n\bar z_k|.
\end{equation*}
Then
\begin{equation*}
\omega (z_n,G_k)\leq \int_{M_0 I_k} \frac{1-|z_n|^2}{|e^{it}-z_n|^2}\; \frac{dt}{2\pi}\leq
C_1^{-2}M_0 \frac{(1-|z_n|^2)(1-|z_k|^2)}{|1-z_n\bar z_k|^2}.
\end{equation*}
It is easy to see from the estimates above that there exists a
universal constant $C_2>0$ such that
\begin{equation*}
\beta (z_n,z_k)\leq C_2 -\log_2
 \frac{(1-|z_n|^2)(1-|z_k|^2)}{|1-z_n\bar z_k|^2},
\end{equation*}
one deduces
\begin{equation*}
\omega (z_n,G_k)\leq  C2^{-\beta (z_n,z_k)}
\end{equation*}
with $C=C_1^{-2} M_0 2^{C_2}$. Hence (\ref{eq3.10}) holds for
points $z_k$ in parts (A) and (B). Therefore
\begin{equation*}
\textrm{(A)} + \textrm{(B)} \leq C\sum_{k\colon \beta(z_k,z_n)\geq
N} 2^{(\eta-1)\beta (z_n,z_k)}.
\end{equation*}
Observe that condition (\ref{1.3}) gives
\begin{equation*}
\sum_{k\colon \beta(z_k,z_n)\leq j} 2^{(\eta-1)\beta
(z_n,z_k)}\leq M 2^{(\eta+\alpha-1)j},
\end{equation*}
for any $j=1,2,\dotsb$ Since $\alpha <1$ one may choose
$0<\eta=\eta(\alpha)<1-\alpha$ so that $\alpha+\eta<1$. So, adding
up for $j\geq N$, one obtains
\begin{equation*}
\textrm{(A)} + \textrm{(B)}\leq CM \frac{2^{(\eta+\alpha-1)N}}{1-2^{\eta+\alpha-1}}.
\end{equation*}
Hence, taking $N>0$ sufficiently large one deduces
\begin{equation*}
\textrm{(A)} + \textrm{(B)} \leq \frac{\delta}{3}.
\end{equation*}

 The estimate of the third term (C) depends on the choice of the constant $\gamma >0$ appearing in the construction of the sets $\{E_n\}$. Fixed $z_n$, consider
\begin{equation*}
U(n)=\left \{ z_k \; \colon \; \beta (z_k,z_n)\geq N, \; z_k\notin
20M_0Q(z_n), \; 2M_0I_k \cap M_0I_n \neq \emptyset \right \}.
\end{equation*}
So $\textrm{(C)} =\sum_{z_k\in U(n)} 2^{\eta \beta (z_k,z_n)}
\omega(z_n, G_k)$.

Observe that if $z_k\in U(n)$,  then  $|z_k|<|z_n|$ and
$z_n\in 3M_0Q(z_k)$. In particular $z_n\in 20M_0 Q(z_k)$ so, by
the construction of the sets $\{G_k\}$, $G_k\subset
M_0I_k\setminus I_n^\gamma (k)$. Hence
\begin{equation*}
\omega(z_n,G_k)  \leq 
\int\limits_{M_0I_k\setminus I_n^\gamma (k)}
 \frac{1-|z_n|^2}{|\xi-z_n|^2} \; \frac{|d\xi|}{2\pi} \leq  
 \int\limits_{\partial  \D\setminus I_n^\gamma (k)}
\frac{1-|z_n|^2}{|\xi-z_n|^2}\;\frac{|d \xi|}{2\pi}
\end{equation*}
and  a change of variable gives an absolute constant $C_3>0$ such
that
\begin{equation}\label{6.6}
\omega(z_n, G_k)\leq
C_3\,(1-|z_n|)\int_{1-|z_n^\gamma(k)|}^{\infty} \frac{dx}{x^2}
\leq C_3\;\frac{1-|z_n|}{1-|z_n^\gamma (k)|}.
\end{equation}
This estimate is worst than (\ref{eq3.10}) which was used for (A) and (B) but it is
good enough for our purposes. The key is that in (C) we sum over ``few"
terms corresponding to the points $z_k\in U(n)$.

Observe that if $z_k\in U(n)$, $z_k$ belongs to the Stolz angle
$\Gamma_n=\Gamma_n(M_0)=\{z\in \D \colon |z-e^{i\argg z_n}|\leq
11M_0(1-|z|)\}$ with vertex $e^{i\operatorname{Arg} z_n}$ and a
certain opening depending on $M_0$. To see this we only need to
observe that $2M_0I_k\cap M_0I_n \neq \emptyset$ implies
$|\operatorname{Arg} z_k- \operatorname{Arg} z_n|\leq
10M_0(1-|z_k|)$ and use this inequality to get
\begin{equation*}
|z_k-e^{i\operatorname{Arg} z_n}|\leq 11 M_0(1-|z_k|).
\end{equation*}

Define $V(n)=\{z_k\in \Gamma_n \colon |z_k|<|z_n|, \beta
(z_k,z_n)\geq N\}$ and then,
\begin{equation*}
\textrm{(C)} =\sum_{z_k\in U(n)} 2^{\eta\beta (z_k,z_n)}
\omega(z_n,G_k)\leq \sum_{z_k\in V(n)} 2^{\eta \beta (z_k,z_n)}
\omega(z_n,G_k).
\end{equation*}
Using inequalities (\ref{6.6}) and (\ref{6.5}) we obtain
\begin{equation*}
\textrm{(C)} \leq C_3\!\!\sum_{z_k\in V(n)}\! \!2^{\eta\beta
(z_k,z_n)} \frac{1-|z_n|}{1-|z_n^\gamma (k)|}\leq C_3
\!\!\sum_{z_k\in V(n)}\!\! 2^{\eta \beta (z_k,z_n)}  \left (
\frac{1-|z_n|}{1-|z_k|}\right )^{C^{-1}\gamma}\!.
\end{equation*}
%
Since $z_n\in 3M_0 Q(z_k)$, Lemma~\ref{lemma35} gives
\begin{equation*}
\left |\beta (z_n,z_k)-\log_2 \frac{1-|z_k|}{1-|z_n|}\right |\leq C (M_0).
\end{equation*}
Hence
\begin{equation*}
 \frac{1-|z_n|}{1-|z_k|}\leq 2^{C (M_0)-\beta (z_n,z_k)}.
\end{equation*}
Therefore
\begin{equation*}
\textrm{(C)}\leq C_3\; 2^{C (M_0)C^{-1}\gamma} \sum_{z_k\in V(n)}
2^{(\eta-C^{-1} \gamma)\beta(z_n,z_k)}.
\end{equation*}
Since the sequence $\{z_n\}$ is separated, there exists
$C_4=C_4(M_0)>0$ such that for any $j\geq 0$, the number of points
$z_k\in V_n$ with $j\leq \beta (z_k,z_n)\leq j+1$ is at most
$C_4$. Hence
\begin{equation*}
\textrm{(C)} \leq C_3\,C_4 2^{C (M_0)C^{-1}\gamma}
\sum_{j=N}^\infty 2^{(\eta-C^{-1} \gamma)j}.
\end{equation*}
Taking $\eta>0$ so small that $\eta-C^{-1}\gamma<0$ and taking $N$ sufficiently large, we deduce
\begin{equation*}
\textrm{(C)} \leq \frac{\delta}{3}.
\end{equation*}
So condition (\ref{3.2}) is satisfied and the proof of Lemma~\ref{lemma33}
 is finished.
\end{proof}

\subsection{Third Step}\label{3.3} On the last step  given a
partition $\{z_n\}=T\cup S$  the sets $\{G_n\}$
 constructed on step 3.2 will be used to find a function
$u=u(T,S)$ satisfying the conditions stated in
Lemma~\ref{lemma31}. This will end the proof of the sufficiency of
condition (\ref{1.3}).

A sequence of points $\{z_n\}$ in the unit disc is called an
interpolating sequence for the space $h^\infty$ of bounded
harmonic functions in the unit disc if for any bounded sequence
$\{w_n\}$ of real numbers there exists $u\in h^\infty$ with
$u(z_n)=w_n$, $n=1,2,\dotsb$. L.~Carleson and J.~Garnett
characterized interpolating sequences for $h^\infty$ as those
sequences $\{z_n\}$ satisfying $\inf_{n\neq m} \beta (z_n,z_m)>0$
and
\begin{equation}\label{3.3}
\sup \frac{1}{\ell (Q)} \sum_{z_n\in Q}1-|z_n|<\infty,
\end{equation}
where the supremum is taken over all Carleson squares of the form
\begin{equation*}
Q=\{ re^{i\theta} \colon 0<1-r<\ell (Q), \quad |\theta-\theta_0|<\ell (Q)\}
\end{equation*}
for some $\theta_0\in [0,2\pi)$. See \cite{CG}, \cite{G1} or
\cite[p. 313]{G2}. 
We next show that a separated sequence $\{z_n\}$ satisfying
(\ref{1.3}) verifies the condition above. Actually it is
sufficient to show (\ref{3.3}) for Carleson squares $Q$ which
contain a point of the sequence  $\{z_n\}$ in its top part
$T(Q)=\{re^{i\theta} \in Q   \colon 1-r > \ell (Q)/2\}$. Let $Q$
be a Carleson square of this type. Let  $z_n\in T(Q)$ and $A(j)=\{
k \colon z_k\in Q ,  j-1\leq \beta (z_k,z_n)< j\}$.
Lemma~\ref{lemma35} tells that for any $k\in A(j)$ the quantity
$1-|z_k|$ is comparable to $2^{-j} \ell (Q)$. Hence condition
(\ref{1.3}) yields
\begin{equation*}
\sum_{k\in A(j)} 1-|z_k|\leq C_12^{-j} \ell (Q)\# A(j)\leq
C_1M2^{(\alpha-1)j} \ell(Q).
\end{equation*}
Since $\alpha<1$, adding up over $j=1,2,\dotsb$, one obtains
(\ref{3.3}). Hence $\{z_n\}$ is an interpolating sequence for
$h^\infty$. Then by the Open Mapping Theorem, there exists a
constant $\gamma =\gamma (\{z_n\})>0$ such that for any partition
of the sequence $\{z_n\}=T\cup S$ there exists $h=h(T,S)\in
h^\infty$ with $\sup \{ |h(z)|\colon z\in \D\} <1$ and
$h(z_n)>\gamma$ for $z_n\in T$ while $h(z_n)<-\gamma$ for $z_n\in
S$. Let $\delta >0$ be a small number to be fixed later and let
$N=N(\delta)$, $\eta=\eta(\delta)$ be the positive constants and
$\{G_n\}$ the pairwise disjoint collection of subsets of the unit
circle given in Lemma~\ref{lemma33}. Let $\varepsilon =\varepsilon
(\delta)$ be a small number to be fixed later which will satisfy
$\varepsilon \delta^{-1} \to 0$ as $\delta$ tends to $0$. Let
$\{w_k\}$ be a sequence of positive numbers satisfying the
compatibility condition (\ref{1.2}). Given a partition
$\{z_n\}=T\cup S$, consider the function $u=u(T,S)\in h^+$ defined
by
\begin{equation*}
u(z)=\sum_k w_k \int_{G_k} P_z (\xi) (1+h(\xi))|d\xi|,
\end{equation*}
where $h=h(T,S)$ and
\begin{equation*}
P_z(\xi)=\frac{1}{2\pi}\; \frac{1-|z|^2}{|\xi-z|^2}
\end{equation*}
is the Poisson kernel. Our goal is to show that $u(z_n)\geq w_n$
for $z_n\in T$ and $u(z_n)\leq w_n$ for $z_n\in S$. For
$n=1,2,\dotsb$, let $A(n)$ be the set of indexes $k$ such that
$\beta (z_k,z_n)\leq N$. Write $u(z_n)=\textrm{(I)} +
\textrm{(II)}$, where
\begin{equation*}
\begin{split}
\textrm{(I)} &= \sum_{k\notin A(n)} \omega_k \int_{G_k} P_{z_n}
(\xi)(1+h(\xi))\,|d \xi|, \\
\textrm{(II)} &= \sum_{k\in A(n)} \omega_k \int_{G_k} P_{z_n}
 (\xi)(1+h(\xi))\,|d \xi|.
\end{split}
\end{equation*}
We first show that
\begin{equation}\label{3.4}
\textrm{(I)}
<2\delta w_n, \quad n=1,2,\dotsb
\end{equation}
Actually if the constant $\varepsilon = \varepsilon ( \delta) >0$
is taken so that $\varepsilon <\eta$, the compatibility condition
(\ref{1.2}) tells that $\textrm{(I)}$ can be bounded by
\begin{equation*}
w_n \sum_{k\notin A(n)} 2^{\eta \beta(z_k,z_n)} 2\omega(z_n,G_k)
\end{equation*}
which, by (\ref{3.2}), is bounded by $2\delta w_n$. Hence
(\ref{3.4}) holds.

For the other term, using that the sets $\{G_n\}$ are pairwise
disjoint and the compatibility condition (\ref{1.2}) we have
\begin{equation*}
\textrm{(II)} =\sum_{k\in A(n)} w_k \int_{G_k} P_{z_n} (\xi)(1+h(\xi))|d\xi|
\leq 2^{\varepsilon N} w_n(1+h(z_n)).
\end{equation*}
Also, since $\sup \{ |h(z_n)|\colon z\in \D\}\leq 1$, the
compatibility condition (\ref{1.2})  and the estimate (\ref{3.1})
yield
\begin{equation*}
\begin{split}
\textrm{(II)}  &\geq w_n 2^{-\varepsilon N} \left (1+h(z_n) -
\int_{\partial \D\setminus \bigcup_{k\in A(n)} G_k} P_{z_n}
(\xi)\,
(1+h(\xi))|d\xi|\right )\\
&\geq 2^{-\varepsilon N} w_n (1+h(z_n)-2\delta).
\end{split}
\end{equation*}
So
\begin{equation*}
2^{-\varepsilon N} w_n (1+h(z_n)-2\delta)\leq \textrm{(II)} \leq
2^{\varepsilon N} w_n (1+h(z_n)).
\end{equation*}
Hence
\begin{enumerate}[(a)]
\item If $z_n\in T$, $h(z_n)\geq \gamma$ and then $u(z_n)\geq
\textrm{(II)} \geq w_n 2^{-\varepsilon N} (1+\gamma -2\delta)$.
\item If $z_n\in S$, $h(z_n)\leq -\gamma$ and then
$u(z_n)=\textrm{(I)}+\textrm{(II)}\leq w_n (2\delta +
2^{\varepsilon N} (1-\gamma))$.
\end{enumerate}
Fixed $\gamma >0$, taking $\delta=\delta(\gamma)>0$ and
$\varepsilon=\varepsilon (\delta, \eta, N)>0$  sufficiently small,
we deduce that $u(z_n)\geq~w_n$ if $z_n\in T$ and $u(z_n)\leq w_n$
if $z_n\in S$. An application of Lemma~\ref{lemma31} concludes the
proof of the sufficiency of condition~(\ref{1.3}).\hspace*{9.5true
cm} $\qed$

\section{Equivalent conditions}\label{sectioncuatro}

In this section several geometric conditions which are equivalent
to (\ref{1.3}) are collected.

\begin{prop}
Let $\{z_n\}$ be a sequence of distinct points in $\D$. Then the
following are equivalent:
\begin{enumerate}[\rm(a)]
\item Condition (\ref{1.3}) holds, that is, there exist constants
$M>0$ and $0<\alpha <1$ such that
\begin{equation*}
\# \{z_j \colon \beta (z_j,z_n)\leq l\} \leq M\,2^{\alpha l}
\end{equation*}
for any $n,l=1,2\dotsb$
\item There exist constants $M_1>0$ and
$0<\alpha <1$ such that
\begin{equation*}
\# \left \{z_j \colon \left |\frac{z_j-z_n}{1-\bar z_nz_j}\right
|\leq r \right \} \leq M_1 (1-r)^{-\alpha},
\end{equation*}
for any $0<r<1$ and any $n=1,2,\dotsb$ \item There exist constants
$M_2>0$ and $0<\alpha<1$ such that
\begin{equation*}
\# \{z_j \in Q(z_n)\colon 2^{-l-1} (1-|z_n|)\leq 1-|z_j|\leq
2^{-l}(1-|z_n|)\} \leq M_2\,2^{\alpha \,l}
\end{equation*}
for any $n,l=1,2, \dotsb$ \item There exist constants $M_3>0$ and
$0<\alpha<1$ such that
\begin{equation*}
\sum_{z_j\in Q(z_n)} (1-|z_j|)^\alpha \leq M_3 (1-|z_n|)^\alpha,
\end{equation*}
for any $n=1,2,\dots$
\end{enumerate}
\end{prop}

\begin{proof}
The equivalence between (a) and (b) follows from the following
obvious observation. Let $z,w\in \D$, then $\beta (z,w)\leq l$ if
and only if
\begin{equation*}
\left |\frac{z-w}{1-\overline w z} \right |= \frac{2^{\beta
(z,w)}-1}{2^{\beta(z,w)}+1}= 1- \frac{2}{2^{\beta(z,w)}+1}\leq
1-\frac{2}{2^l+1}
\end{equation*}
Assume (a) holds. Fix two positive integers $n,l$. Let $z_j\in
Q(z_n)$ satisfying
\begin{equation*}
2^{-l-1} (1-|z_n|)\leq 1-|z_j|\leq 2^{-l}(1-|z_n|).
\end{equation*}
Applying Lemma \ref{lemma35} one shows that there exists a
universal constant $C>0$ such that
\begin{equation*}
|\beta (z_n,z_j)-l|\leq C.
\end{equation*}
Hence
\begin{multline*}
\left \{z_j\in Q (z_n) \colon 2^{-l-1}  (1-|z_n|)\leq 1-|z_j| \leq
2^{-l} (1-|z_n|)\right \}\\ \subseteq \{z_j\colon \beta
(z_j,z_n)\leq l+C\}
\end{multline*}
and condition (\ref{1.3}) gives (c). Adding up over $l=1,2,\dots$
one shows that (c) implies (d). Assume (d) holds and let us show
condition (\ref{1.3}). By conformal invariance one may assume
$z_n=0$. So condition (d) tells
\begin{equation*}
\sum_{j=1}^\infty (1-|z_j|)^\alpha \leq M_3.
\end{equation*}
Since $\beta (z_j,0)\leq l$ implies
\begin{equation*}
1-|z_j|\geq \frac{2}{2^{l}+1},
\end{equation*}
we deduce
\begin{equation*}
\# \{z_j \colon \beta (z_j,0)\leq l\} \leq M_3\, \left
(\frac{2}{2^l+1} \right )^{-\alpha}
\end{equation*}
which gives (\ref{1.3})
\end{proof}

As mentioned in the introduction, condition (\ref{1.3}) tells how
dense is the sequence when one looks at it from a point of the
sequence. It is worth mentioning that one can not take as a base
point an arbitrary point in the unit disc. This follows from the
following example of two separated interpolating sequences for
$h^+$ which will be called $Z_1$, $Z_2$ so that $\inf \{\beta
(z,\xi)\colon z\in Z_1,\; \xi\in Z_2\}>0$ but such that the union
$Z_1\cup Z_2$ is not an interpolating sequence for $h^+$. For
instance one may take $Z_1=\{r_k\}$ where $r_1=0$, $r_k \to 1$ and
$\beta(r_k,r_{k+1})\to \infty$ as $k\to\infty$. For each
$k=1,2,\dotsb$, choose points $\{z_1^{(k)}, \dotsb,
z_{N(k)}^{(k)}\}$, $N(k)=2^{n_k}$, equally distributed in the
hyperbolic cercle centered at $r_k$ of hyperbolic radius $n_k$.
Here $n_k \to \infty$ as $k\to \infty$ in such a way that
$n_k<\beta (r_k, r_{k+1})/4$. Let $Z_2=\{ z_i^{(k)}\colon
i=1,\dotsb, N(k), k=1,2,\dotsb\}$. It can be shown that $Z_1$ and
$Z_2$ satisfy condition (\ref{1.3}) with the exponent
$\alpha=1/2$, while $Z_1\cup Z_2$ does not fulfill (\ref{1.3}) for
any $0<\alpha <1$ because the number of points in $Z_2$ at
hyperbolic distance $n_k$ from the point $r_k\in Z_1$ is
$2^{n_k}$.

\section{An interpolation problem for bounded
Analytic Functions without zeros}\label{sectioncinco}

Let $\Hh^\infty$ denote the algebra of bounded analytic functions
in the unit disc $\D$. Let $(\Hh^\infty)^*$ be the subalgebra of
$\Hh^\infty$ which consists on the functions in $\Hh^\infty$
without zeros in $\D$. If $f\in (\Hh^\infty)^*$ then $\log \left (
\|f\|_\infty / |f(z)|\right )\in h^+$. So if $\{z_n\}$ is a
sequence in $\D$ and $t_n=\log \left ( \|f\|_\infty
/|f(z_n)|\,\right )$, Harnack's inequality tells that
\begin{equation*}
|\log t_n -\log t_m|\leq \beta (z_n,z_m), \qquad n,m=1,2,\dotsb
\end{equation*}
So, as before, we may consider a notion of interpolating sequence.

\begin{defn}\label{def4.1}
A sequence of points $\{z_n\}$ in the unit disc is called an interpolating sequence for $(\Hh^\infty)^*$ if there exist
 constants $\varepsilon >0$ and $0<C<\infty$ such that for any sequence
  of non-vanishing complex values $\{w_n\}$, $|w_n| < C, n=1,2,\dotsb $, satisfying
\begin{equation}\label{4.1}
\left | \log \left ( \log \left ( \frac{C}{|w_n|}\right )\right )
-\log \left (\log \left (\frac{C}{|w_m|}\right )\right )\right |
\leq \varepsilon \beta (z_n,z_m),\; n,m=1,2,\dotsb
\end{equation}
there exists a function $f\in (\Hh^\infty)^*$ with $f(z_n)=w_n$,
$n=1,2,\dotsb$
\end{defn}

The characterization of the interpolating sequences for
$(\Hh^\infty)^*$ is given in the following result.

\begin{thm}\label{thmbaf}
A separated sequence $\{z_n\}$ of points in the unit disc is
interpolating for $(\Hh^\infty)^*$  if and only if there
exist constants $M>0$ and $0<\alpha <1$ such that
\begin{equation}\label{4.2}
\# \{ z_j \colon \beta (z_j,z_n)\leq \ell \}\leq  M2^{\alpha \ell}
\textrm{ for any } n,\ell=1,2,\dotsb
\end{equation}
\end{thm}

 \begin{proof}[Proof of Theorem~\ref{thmbaf}]
Let us start by showing the neccessity of condition (\ref{4.2}).
Given a separated interpolating sequence  $\{z_n\}$ for
$(\Hh^\infty)^*$  consider the constants $\varepsilon
>0$ and $C<\infty$ given in definition~\ref{def4.1}.
Define the sequence of positive values $t_n=2^{\varepsilon\beta
(0,z_n)}$, $\; n=1,2,\dotsb$ It is clear that
\begin{equation*}
\left | \log_2 t_n -\log_2 t_m\right | \leq \varepsilon\beta
(z_n,z_m), \quad n,m=1,2,\dotsb.
\end{equation*}
Then, if we consider a sequence of complex values $\{w_n\}$ with
$t_n=\log\left ( C/|w_n|\right )$,  we  have $\sup_n |w_n|\leq C$
 and furthermore
$\{w_n\}$ satisfies condition~(\ref{4.1}). So, there exists a
function $f\in \Hh^\infty$ without zeros with $f(z_n)=w_n$,
$n=1,2,\dotsb$ The function $v(z)=\log \left
(\frac{C}{|f(z)|}\right )$ is a harmonic function, $v(z)\geq \log
(C/\|f\|_\infty):= -k_1$, and interpolates the values $\{t_n\}$ at
the points $\{z_n\}$. So, $u(z)=v(z)+k_1\in h^+(\D)$ and
$u(z_n)=t_n+k_1=2^{\varepsilon \beta (0,z_n)}+k_1$, $n=1,2,,
\dotsb$ Now, arguing as in the proof of the necessity of
Theorem~\ref{thm1}, we can conclude that there exist constants
$M>0$ and $0<\alpha <1$ such that
\begin{equation*}
\#\{z_j \colon \beta (z_j,z_n)\leq \ell\} \leq M 2^{\alpha \ell} \textrm{ for any } n,\ell=1,2,\dotsb
\end{equation*}

Let us now show the sufficiency of condition (\ref{4.2}). Given a
separated sequence $\{z_n\}$ satisfying (\ref{4.2}) and $\{w_n\}$
satisfying (\ref{4.1}) for some $\varepsilon, C$, consider
$t_n=\log \frac{C}{|w_n|}$. We can take $C>\|w_n\|_\infty$. Then
obviously $\{t_n\}$ satisfies the compatibility condition
(\ref{1.2}). So, there exists a function $u\in h^+ (\D)$ with
$u(z_n)=\log \frac{C}{|w_n|}$, for $n=1,2,\dotsb$ Consider
$u_0(z)=u(z)-\log (C)$ and let $\widetilde{u_0}(z)$ be the
 harmonic conjugate function of $u_0(z)$. Then
$e^{-(u_0+i\widetilde{u_0})}$ is a bounded analytic function that
interpolates de values $\{|w_n|\gamma_n\}$ at the points
$\{z_n\}$, where $\gamma _n=e^{-i\widetilde u_0(z_n)}$,
$n=1,2,\dotsb$. The sequence $\{z_n\}$ is separated and satisfies
condition (\ref{1.3}), so it is an interpolating sequence for
$\Hh^\infty$ (see \cite{C1} or \cite{G2}). So there exists a
bounded analytic function $g(z)$ such that
$g(z_n)=-\operatorname{Arg} (\gamma_n)+\operatorname{Arg} (w_n)$
and then the function $h(z)=e^{-u_0-i\widetilde{u_0}}e^{ig}$ is a
bounded analytic function without zeros with $h(z_n)=w_n$ for any
$n=1,2,\dotsb$
\end{proof}

\section{Higher Dimensions}\label{sectionseis}

Let $h^\infty (\Rr_+^{d+1})$ be the space of bounded harmonic
functions in the upper-half space $\Rr_+^{d+1}=\{(x,y)\colon x\in
\Rr^d, y>0\}$. A sequence of points $\{z_n\}\subset \Rr_+^{d+1}$
is called an interpolating sequence for $h^\infty (\Rr_+^{d+1})$
if for any bounded sequence $\{w_n\}$ of real numbers there exists
$u\in h^\infty (\Rr_+^{d+1})$ with $u(z_n)=w_n$, $n=1,2,\dotsb$.
When the dimension $d>1$, there is no complete geometric
description of the interpolating sequences for $h^\infty
(\Rr_+^{d+1})$. In \cite{C1} and \cite{CG}, L.~Carleson and
J.~Garnett proved the following result.

\begin{thm}\cite{C1}, \cite{CG} Let $\{z_n=(x_n,y_n)\}$ be a sequence of points in $\Rr_+^{d+1}$, $d>1$.
\begin{enumerate}[\rm(a)]
\item Assume $\{z_n\}$ is an interpolating sequence for $h^\infty
(\Rr_+^{d+1})$. Then
\begin{equation}\label{5.1}
\inf_{n\neq m} \beta (z_n,z_m)>0
\end{equation}
and there exists a constant $C=C(\{z_n\})$ such that
\begin{equation}\label{5.2}
\sum_{z_n\in Q} y_n^d \leq C \ell (Q)^d
\end{equation}
for any Carleson cube $Q$ of the form
\begin{equation*}
Q=\{(x,y)\in \Rr_+^{d+1} \colon |x-x_0|<\ell (Q), \quad
0<y<\ell(Q)\},
\end{equation*}
where $x_o\in \Rr^d$.
\item Assume $\{z_n\}$ satisfies the two
conditions (\ref{5.1}) and (\ref{5.2}) above. Then $\{z_n\}$ can
be splitted into a finite number of disjoint subsequences
$\Lambda_j$, $j=1,\dotsb, N$, that is,
\begin{equation*}
\{z_n\} =\Lambda_1\cup \dotsb \cup \Lambda_N,
\end{equation*}
such that $\Lambda_i\cup \Lambda_j$ is an interpolating sequence for $h^\infty (\Rr_+^{d+1})$ for any $i,j=1,\dotsb, N$.
\end{enumerate}
\end{thm}

Here $\beta (z,w)$ denotes the hyperbolic distance between the points $z,w\in \Rr_+^{d+1}$,
\begin{equation*}
\beta (z,w)=\log_2 \frac{1+\rho (z,w)}{1-\rho (z, w)},
\end{equation*}
where $\rho (z,w)= |z-w| / |z-\bar w|$ and $\bar w=(w_1, \dotsb,
w_d,-w_{d+1}).$

Moreover in \cite{CG}, the authors present several geometric
conditions on the sequence $\{z_n\}$ which imply that $\{z_n\}$ is
an interpolating sequence for $h^\infty (\Rr_+^{d+1})$. However it
is not known if the two necessary conditions (\ref{5.1}) and
(\ref{5.2}) are sufficient.  Related interpolation problems have
been considered in \cite{A} and \cite{D}. The situation for
interpolating sequences for the space $h^+(\Rr_+^{d+1})$ of
positive harmonic functions in $\Rr_+^{d+1}$ is quite similar. A
sequence of points $\{z_n\}\subset \Rr_+^{d+1}$ will be called an
interpolating sequence for $h^+(\Rr_+^{d+1})$ if there exists a
constant
 $\varepsilon =\varepsilon (\{z_n\}) >0$ such that for any sequence
 $\{w_n\}$ of positive values satisfying
\begin{equation*}
|\log_2 w_n -\log_2 w_m|\leq \varepsilon \beta (z_n,z_m), \quad
n,m=1,2,\dotsb,
\end{equation*}
there exists a function $u\in h^+(\Rr_+^{d+1})$ with $u(z_n)=w_n$,
$n=1,2, \dotsb$

As before, a sequence of points $\{z_n\}\subset \Rr_+^{d+1}$ is
called separated if $\inf_{n\neq m} \beta (z_n,z_m)>0$.

\begin{thm}
Let $\{z_n\}$ be a separated sequence of points in the upper-half space
$\Rr_+^{d+1}$, $d>1$.
\begin{enumerate}[\rm(a)]
\item Assume that $\{z_n\}$ is an interpolating sequence for
$h^+ (\Rr_+^{d+1})$. Then there exist constants $M>0$ and $0<\alpha <1$
such that
\begin{equation}\label{5.3}
\# \{z_k \colon \beta (z_k,z_n)\leq l \}\leq M2^{\alpha d l},
\quad l, n=1,2,\dotsb
\end{equation}
 \item Assume that $\{z_n\}$ satisfies the condition (\ref{5.3}) above.
 Then $\{z_n\}$ can be splitted into a finite number of disjoint
 subsequences $\Lambda_i$, $i=1, \dotsb,
 N$,
\begin{equation*}\label{45}
\{z_n\}=\Lambda_1 \cup \dotsb \cup \Lambda_n,
\end{equation*}
such that $\Lambda_i\cup \Lambda_j$ is an interpolating sequence for
$h^+ (\Rr_+^{d+1})$ for any $i,j=1,\dotsb, N$
\end{enumerate}
\end{thm}

The proof of (a) follows the same lines of the proof of the
necessity in Theorem~\ref{thm1}. The first two steps \ref{3.1} and
\ref{3.2} of the proof of the sufficiency in Theorem~\ref{thm1}
can be extended to several variables. However the third step 3.3
can not be fulfilled because we have not been able to show that a
separated sequence satisfying condition (\ref{5.3}) is an
interpolating sequence for $h^\infty (\Rr_+^{d+1})$. Since it is
clear that (\ref{5.3}) implies (\ref{5.2}), applying the result of
L.\ Carleson and J.\ Garnett, the sequence $\{z_n\}$ can be
splitted into a finite number of disjoint subsequences
$\Lambda_1,\dotsb,\Lambda_N$ such that $\Lambda_i \cup \Lambda_j$
is an interpolating sequence for $h^\infty (\Rr_+^{d+1})$,
$i,j=1,\dotsb, N$. Arguing as in step 3.3 of the proof of the
suffiency, one can show that for any $i,j=1,\dotsb,N$,
 the sequence $\Lambda_i \cup \Lambda_j$ is an interpolating sequence
 for $h^+ (\Rr_+^{d+1})$.

It is worth mentioning that we have not been able to prove that a
separated sequence verifying (\ref{5.3}) is interpolating for
$h^+(\Rr_+^{d+1})$, when $d>1$.



\begin{thebibliography}{AAA}
\bibitem[A]{A} E.\ Amar, \emph{Suites d'interpolation harmoniques},
 J.\ Anal. Math., \textbf{32} (1977), 197--211.
%
\bibitem[BN]{BN} B.\ Boe \& A.\ Nicolau, \emph{Interpolation by functions
in the Bloch space}, J.\ Anal.\ Math., \textbf{94} (2004),
171--194.
%
\bibitem[C1]{C1} L.\ Carleson, \emph{An interpolation problem for bounded
analytic functions}, Amer.\ J.\ Math.,  {\bf 80} (1958), 921--930.
%
\bibitem[C2]{C2} L.\ Carleson, \emph{A moment problem and harmonic
interpolation}, preprint, Institut Mittag-Leffler, 1972
%
\bibitem[CG]{CG} L.\ Carleson \& J.\ Garnett, \emph{Interpolating sequences
and separation properties}, J.\ Anal.\ Math., {\bf 28} (1975),
273--299.
%
\bibitem[D]{D} K.\ Dyakonov, \emph{Moment problems for bounded functions},
Comm.\ Anal.\ Geom.,  {\bf 2(4)} (1994), 533--562.
%
\bibitem[DN]{DN} K.\ Dyakonov \& A.\ Nicolau, \emph{Free interpolation by
non-vanishing analytic functions}, to appear in Trans.\ Amer.\ Math.\ Soc.
%
\bibitem[G1]{G1} J.\ Garnett, \emph{Interpolating sequences for
bounded harmonic functions}, Indiana Univ. Math.\ J., \textbf{21}
(1971/1972), 187--192.
%
\bibitem[G2]{G2} J.B.\ Garnett, \emph{Bounded Analytic Functions}, Pure and
Applied Mathematics, 96. Academic Press, Inc., New York-London,
1981.
%
\bibitem[H]{H} Hall, \emph{Sur la mesure harmonique de certains ensembles},
Arkiv f\"or Matematik, Astronomi och Fysik,  {\bf 25A(28)} (1937).
%
\bibitem[HL]{HL} J.B.\ Hiriart-Urruty \& C.\ Lemar\'echal,
\emph{Convex Analysis and Minimization Algorithms I.
Fundamentals}, Grundlehren der Mathematischen Wissenschaften 305,
Springer-Verlag, Berlin, Heidelberg, New York, 1993.
%
\bibitem[MS]{MS} D.E.\ Marshall \& C.\ Sundberg, \emph{Harmonic measure and radial
projection}, Trans.\ Amer.\ Math.\ Soc.,  {\bf 316(1)} (1989),
81--95.
%
Graduate Texts in Mathematics \textbf{149} Springer-Verlag, New
York, 1994.
%
\bibitem[S]{S} K.\ Seip, \emph{Interpolating and sampling in spaces of
analytic functions}, University Lectures Series, \textbf{33}
American Mathematical Society, Providence, RI, 2004. 


\end{thebibliography}
\end{document}